\documentclass[11 pt]{article}
\usepackage[english]{babel}
\usepackage[T1]{fontenc}
\usepackage{pifont,setspace,hyperref,amsthm,amsmath,amssymb,color,multirow,graphicx,subfigure, color }
\usepackage{enumitem}
\usepackage{inputenc}
\usepackage{amsmath,latexsym,amsthm}
\usepackage{graphicx}
\usepackage{amssymb}
\begin{document}
\begin{center}
\textbf{\Large \textsc {Tricomi type boundary value problem with integral conjugation condition for a mixed type equation with Hilfer fractional operator}}\\[0.2 cm]
\textbf{Karimov E.~T.~\footnote{V.I.Romanovskiy Institute of Mathematics of Uzbekistan Academy of Sciences, Tashkent, Uzbekistan. E-mail: erkinjon@gmail.com}
\footnote{Ghent University, Ghent, Belgium. E-mail: ekarimov@ugent.be}
\footnote{Sultan Qaboos University, FracDiff Research Group, Muscat, Oman. E-mail: erkinjon@squ.edu.om}}
\end{center}
\begin{abstract}
In this work, The Tricomi type boundary problem with integral conjugation condition on the type-changing line for the mixed type equation with Hilfer fractional differential operator has been considered. Using method of integral equations, energy integral's method, a unique solvability of the formulated problem has been proved.
\end{abstract}
\vspace{0.5 cm}


\section{Introduction}

Fractional Calculus is developing intensively due to both practical \cite{uch} and theoretical importance \cite{kilb}. Fractional analogues of essential equations such as diffusion, wave were studied involving different fractional differential and integral operators\cite{gorluch},\cite{sandev},\cite{pskhu},\cite{salti}.  We will omit a huge amount of works devoted to the studying of direct and inverse problems for partial differential equations (PDEs) with fractional order operators and note only studies closely related to the present topic. Several boundary value problems (BVPs) for mixed type equations with Riemann-Liouville fractional differential operator (FDO) were studied for unique solvability in works \cite{kilbrep},\cite{gekk},\cite{kad}. BVP with integral form conjugation conditions for PDEs with both Riemann-Liouville and Caputo FDOs were subject of series of investi\-ga\-tions \cite{akhatkar}, \cite{cabkar}, \cite{berkar}, \cite{agarkar}, \cite{abd}. In these works, authors used an explicit solution of certain BVP for fractional diffusion equation studied by A.Pskhu \cite{pskhu}.

For the first time, generalized Riemann-Liouville FDO which is also called as Hilfer FDO introduced by Hilfer \cite{hilf}. The Cauchy and some BVPs for ODEs and PDEs with Hilfer FDO investigated by many authors, for instance \cite{luch}, \cite{sandev}, \cite{fur}, \cite{malik}.

In this paper, we are aimed to study BVP with integral form conjugation condition in a mixed domain consisted of characteristic triangle and rectangle, for a mixed type PDE with diffusion equation involving Hilfer FDO.
\section{Formulation of a problem and main functional relations}
\subsection{Formulation of a problem}

Let us consider the following mixed type equation
\begin{equation}\label{eq1}
0=\left\{\begin{array}{l} {k}u_{xx}-D_{0t}^{\alpha,\mu}u,\,\,t>0\\ u_{xx}-u_{tt},\,\,t<0\end{array}\right.
\end{equation}
in a mixed domain $\Omega=\Omega_1\cup\Omega_2\cup AB$. Here ${k}=const>0$, $\Omega_1=\left\{(x,t):\, 0<x<l,\,0<t<T\right\}$, $\Omega_2=\left\{(x,t):\,-t<x<t+l,\,-l/2<t<0\right\}$ , $AB=\left\{(x,t):\,0<x<l,\,t=0\right\}$,
\[
D_{0t}^{\alpha,\mu}f=I_{0t}^{\mu(1-\alpha)}\frac{d}{dt}I_{0t}^{(1-\mu)(1-\alpha)}f(t)
\]
is the Hilfer FDO of order $\alpha\, (0<\alpha\leq 1)$ and of type $\mu\, (0\leq \mu\leq 1)$ \cite{hilf}, where
\[
I_{0t}^{\alpha}f(t)=\frac{1}{\Gamma(\alpha)}\int\limits_0^t\dfrac{f(z)dz}{(t-z)^{1-\alpha}}
\]
is the Riemann-Liouville fractional integral of order $\alpha\, (\Re (\alpha)>0)$ \cite{kilb}.
\begin{figure}[h]
\centering
\includegraphics[width=0.4\linewidth]{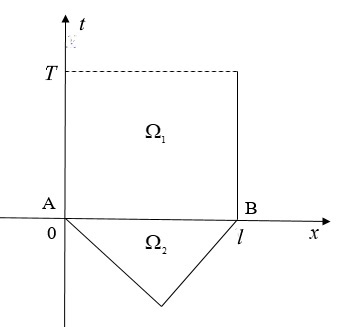}
\end{figure}

Tricomi type BVP for Eq.(\ref{eq1}) in $\Omega$ can be formulated as follows:

\textbf{Problem.} To find a function $u(x,t)$ which is continuous in $\overline{\Omega}\backslash AB$, its Hilfer derivative is continuous in $\Omega_1$ and it has continuous second order partial derivatives in $\Omega_2$, and it satisfies Eq. (\ref{eq1}) in $\Omega$ together with boundary conditions
\begin{equation}\label{eq2}
u(0,t)=0,\,\,u(l,t)=0,\,0\leq t\leq T,
\end{equation}
\begin{equation}\label{eq3}
u(x/2,-x/2)=\psi(x),\,\,\,0\leq x\leq l,
\end{equation}
conjugation conditions on $AB$
\begin{equation}\label{eq4}
\lim\limits_{t\rightarrow +0} t^{(1-\mu)(1-\alpha)}u(x,t)=u(x,-0),\,\,0\leq x\leq l,
\end{equation}
\begin{equation}\label{eq5}
\begin{array}{l}
\lim\limits_{t\rightarrow +0} t^{1-\alpha}\left(t^{(1-\mu)(1-\alpha)}u(x,t)\right)_t=\gamma_1u_t(x,-0)+\gamma_2\int\limits_0^xu_t(z,-0)P(x,z)dz+\\
\gamma_3\int\limits_x^lu_t(z,-0)Q(x,z)dz,\,\,0< x< l.
\end{array}
\end{equation}
Here $\gamma_i\, (i=\overline{1,3})$ are real constants, $\psi(x)$, $P(\cdot,\cdot), Q(\cdot,\cdot)$ are given functions such that $\psi(0)=0$, $\sum\limits_{i=1}^3\gamma_i^2>0$.

\subsection{Main functional relations}

Let us introduce a notation
\begin{equation}\label{eq6}
\tau_1(x)=\lim\limits_{t\rightarrow+0}t^{(1-\mu)(1-\alpha)}u(x,t),\,\,0\leq x\leq l.
\end{equation}

Solution of the Eq.(\ref{eq1}) in $\Omega_1$ which satisfies conditions (\ref{eq2}), (\ref{eq6}) can be written as follows \cite{sandev}:
\begin{equation}\label{eq7}
u(x,t)=\frac{2}{l}\Gamma[\alpha+\mu(1-\alpha)]\int\limits_0^l\tau_1(\xi)t^{(1-\mu)(\alpha-1)}\sum\limits_{n=1}^\infty E_{\alpha,\alpha+\mu(1-\alpha)}\left(-k\lambda_n^2t^\alpha\right)\sin(\lambda_nx)\sin(\lambda_n\xi)d\xi,
\end{equation}
where $\lambda_n=\frac{n\pi}{l}$,
\begin{equation}\label{eq8}
E_{\alpha,\beta}(z)=\sum\limits_{j=0}^\infty\frac{z^j}{\Gamma(\alpha j+\beta)}
\end{equation}
is a two-parameter Mittag-Leffler function \cite{Podl}.

Using representations (\ref{eq7}) and (\ref{eq8}), we evaluate $t^{1-\alpha}\left(t^{(1-\mu)(1-\alpha)}u(x,t)\right)_t$:
\[
t^{1-\alpha}\left(\frac{2}{l}\Gamma[\alpha+\mu(1-\alpha)]\int\limits_0^l\tau_1(\xi)\sum\limits_{n=1}^\infty E_{\alpha,\alpha+\mu(1-\alpha)}\left(-k\lambda_n^2t^\alpha\right)\sin(\lambda_nx)\sin(\lambda_n\xi)d\xi\right)_t=
\]
\[
\frac{2}{l}t^{1-\alpha}\Gamma[\alpha+\mu(1-\alpha)]\int\limits_0^l\tau_1(\xi)\sum\limits_{n=1}^\infty \frac{d}{dt}\left(\frac{1}{\Gamma(\alpha+\mu(1-\alpha))}+\frac{-k\lambda_n^2t^\alpha}{\Gamma(2\alpha_\mu(1-\alpha))}+\frac{(-k\lambda_n^2t^\alpha)^2}{\Gamma(3\alpha_\mu(1-\alpha))}+...\right)\times
\]
\[
\sin(\lambda_nx)\sin(\lambda_n\xi)d\xi=-k \alpha \Gamma[\alpha+\mu(1-\alpha)]\frac{2}{l}\int\limits_0^l\tau_1(\xi)\sum\limits_{n=1}^\infty \lambda_n^2\left(\frac{1}{\Gamma(2\alpha+\mu(1-\alpha))}+\frac{-2k\lambda_n^2t^\alpha}{\Gamma(3\alpha+\mu(1-\alpha))}+\right.
\]
\[
\left.\frac{3(-k\lambda_n^2t^\alpha)^3}{\Gamma(4\alpha+\mu(1-\alpha))}+...\right)\sin(\lambda_nx)\sin(\lambda_n\xi)d\xi
\]

We introduce another notation, namely
\begin{equation}\label{eq9}
\nu_1(x)=\lim\limits_{t\rightarrow+0}t^{1-\alpha}\left(t^{(1-\mu)(1-\alpha)}u(x,t)\right)_t,\,\,0<x<l.
\end{equation}

Considering above-given evaluations from (7) we obtain the following functional relation on $AB$ deduced from $\Omega_1$ as $t\rightarrow +0$:
\begin{equation}\label{eq10}
\nu_1(x)=\frac{k\alpha\Gamma[\alpha+\mu(1-\alpha)]}{\Gamma[2\alpha+\mu(1-\alpha)]}\tau_1''(x),\,\,\,0<x<l.
\end{equation}

Now we will establish another functional relation on $AB$ which will be reduced from $\Omega_2$. For this aim, we use a solution of the following Cauchy problem:
\[
u_{xx}-u_{tt}=0,\,\,u(x,-0)=\tau_2(x),\,0\leq x\leq l,\,\,u_t(x,-0)=\nu_2(x),\,\,0<x<l,
\]
which has a form \cite{Evans}
\begin{equation}\label{eq11}
u(x,t)=\frac{1}{2}\left[\tau_2(x+t)+\tau_2(x-t)+\int\limits_{x-t}^{x+t}\nu_2(z)dz\right].
\end{equation}
We substitute (\ref{eq11}) into (\ref{eq3}) and deduce
\begin{equation}\label{eq12}
\nu_2(x)=\tau_2'(x)-2\psi'(x),\,\,0<x<l.
\end{equation}
\section{Existence and uniqueness results}
\subsection{Existence of a solution}

Considering conjugation conditions (\ref{eq4}), (\ref{eq5}), from functional relations (\ref{eq10}) and (\ref{eq12}) we get
\begin{equation}\label{eq13}
\tau_1''(x)-A\tau_1'(x)=F_1(x),\,\,0<x<l,
\end{equation}
where $A=\gamma_1\Gamma[2\alpha+\mu(1-\alpha)]/(k\alpha\Gamma[\alpha+\mu(1-\alpha)])$.

Boundary conditions (\ref{eq2}) yield
\begin{equation}\label{eq14}
\tau_1(0)=0,\,\,\,\tau_1(l)=0.
\end{equation}
Solution of (\ref{eq13})-(\ref{eq14}) can be written as \cite{dzhur}
\begin{equation}\label{eq15}
\tau_1(x)=\int\limits_0^lG_0(x,\xi)F_1(\xi)d\xi,
\end{equation}
where
\begin{equation}\label{eq16}
\begin{array}{l}
F_1(x)=\frac{\Gamma[2\alpha+\mu(1-\alpha)]}{k\alpha\Gamma[\alpha+\mu(1-\alpha)]}\left\{\gamma_2\int\limits_0^x\tau_1'(z)P(x,z)dz+\gamma_3\int\limits_0^x\tau_1'(z)Q(x,z)dz-\right.\\
\left.2\gamma_1\psi'(x)-2\gamma_2\int\limits_0^x\psi'(z)P(x,z)dz-2\gamma_2\int\limits_x^l\psi'(z)Q(x,z)dz \right\},
\end{array}
\end{equation}
\[
G_0(x,\xi)=\frac{1}{A\left[e^{Ax}-e^{A(x-l)}\right]}\left\{
\begin{array}{l}
\left(1-e^{A\xi}\right)\left(1-e^{A(x-l)}\right),\,\,0\leq \xi\leq x,\\
\left(1-e^{A(\xi-l}\right)\left(1-e^{Ax}\right),\,\,x\leq \xi\leq l.
\end{array}
\right.
\]
Substituting (\ref{eq16}) into (\ref{eq15}), after integration by parts, we will get
\begin{equation}\label{eq17}
\tau_1(x)-\int\limits_0^l\tau_1(\xi)K(x,\xi)d\xi=F_2(x),\,\,0\leq x\leq l,
\end{equation}
where
\[
\begin{array}{l}
K(x,\xi)=\frac{\Gamma[2\alpha+\mu(1-\alpha)]}{k\alpha\Gamma[\alpha+\mu(1-\alpha)]}\times\\
\left\{G_0(x,\xi)[\gamma_2P(\xi,\xi)-\gamma_3Q(\xi,\xi)]-
\int\limits_\xi^l\left[\gamma_2\frac{\partial P(z,\xi)}{\partial\xi}+\gamma_3\frac{\partial Q(z,\xi)}{\partial\xi}\right]G_(x,z)dz,\right\}
\end{array}
\]
\[
F_2(x)=\frac{-2\Gamma[2\alpha+\mu(1-\alpha)]}{k\alpha\Gamma[\alpha+\mu(1-\alpha)]}\int\limits_0^lG_0(x,\xi)\left[\gamma_1\psi(\xi)+\gamma_2\int\limits_0^\xi\psi'(z)P(\xi,z)dz+\gamma_3\int\limits_\xi^l\psi'(z)Q(\xi,z)dz\right]d\xi.
\]
If $K(x,\xi)$ is continuous or has a weak singularity and $F_2(x)$ is continuously differentiable, then a solution of the second kind Fredholm integral equation (\ref{eq17}) can be represented via resolvent-kernel \cite{Krasn}:
\begin{equation}\label{eq18}
\tau_1(x)=F_2(x)+\int\limits_0^l F_2(\xi)R(x,\xi)d\xi,
\end{equation}
where $R(x,\xi)$ is a resolvent-kernel of $K(x,\xi)$.

We have reduced considered problem to the second kind Fredholm integral equation, which is solvable under the certain conditions to the given data, but it might be not unique. Hence, we have to prove a uniqueness of the formulated problem separately.

\subsection{A uniqueness of the solution}
We multiply equality (\ref{eq10}) by $\tau_1(x)$ and integrate along $AB$:
\[
\int\limits_0^l \tau_1(x)\tau_1''(x)dx-\frac{\Gamma[2\alpha+\mu(1-\alpha)]}{k\alpha\Gamma[\alpha+\mu(1-\alpha)]}\int\limits_0^l\tau_1(x)\nu_1(x)dx=0.
\]
Considering $\int\limits_0^l \tau_1(x)\tau_1''(x)dx=-\int\limits_0^l\left(\tau_1'(x)\right)^2dx$, we deduce
\begin{equation}\label{eq19}
\frac{\Gamma[2\alpha+\mu(1-\alpha)]}{k\alpha\Gamma[\alpha+\mu(1-\alpha)]}\int\limits_0^l\tau_1(x)\nu_1(x)dx+\int\limits_0^l\left(\tau_1'(x)\right)^2dx=0.
\end{equation}

Let us first consider the following integral
\begin{equation}\label{eq20}
\mathbb{I}=\int\limits_0^l\tau_1(x)\nu_1(x)dx,
\end{equation}
where $\tau_1(x)$ and $\nu_1(x)$ are defined by (\ref{eq6}) and (\ref{eq9}), respectively.

Considering (\ref{eq5}) and (\ref{eq12}) at $\psi(x)\equiv 0$, after integration by parts we get
\begin{equation}\label{eq21}
\mathbb{I}=\int\limits_0^l\tau_1^2(x)\left[\gamma_2P(x,x)+\gamma_3Q(x,x)\right]dx-\int\limits_0^l\tau_1(x)\left[\gamma_2\int\limits_0^x\tau_1(z)\frac{\partial P(x,z)}{\partial z}-\gamma_2\int\limits_x^l\tau_1(z)\frac{\partial Q(x,z)}{\partial z}\right]dx.
\end{equation}

Suppose that
\begin{equation}\label{eq22}
\frac{\partial P(x,z)}{\partial z}=-P_1(x)P_1(z),\,\frac{\partial Q(x,z)}{\partial z}=-Q_1(x)Q_1(z),
\end{equation}
then from (\ref{eq21}) it follows that
\begin{equation}\label{eq23}
\mathbb{I}=\int\limits_0^l\tau_1^2(x)\left[\gamma_2P(x,x)+\gamma_3Q(x,x)\right]dx+\frac{\gamma_2\Phi_1^2(l)}{2}+\frac{\gamma_3\Phi_2^2(0)}{2},
\end{equation}
where
\[
\Phi_1(x)=\int\limits_0^x\tau_1(z)P_1(z)dz,\,\,\,\,\Phi_2(x)=\int\limits_x^l\tau_1(z)Q_1(z)dz.
\]
If we suppose that
\begin{equation}\label{eq24}
\gamma_2\geq0,\,\gamma_3\geq0,\,P(x,x)\geq0,\,Q(x,x)\geq0,
\end{equation}
from (\ref{eq23}) we will get $\mathbb{I}\geq0$.

Since $\Gamma(z)>0$ for all $z>0$, then for $k>0$ from (\ref{eq19}) we will have $\tau_1(x)\equiv 0$. Further, considering solution of the first BVP for Eq.(\ref{eq1}) in $\Omega_1$ \cite{sandev}, we will get $u(x,t)\equiv0$ in $\overline{\Omega_1}$. Due to (\ref{eq4}), one can easily deduce that $u(x,t)\equiv0$ in $\overline{\Omega}$.

Finally, we are able now formulate our result as the following

\textbf{Theorem.} If $\psi(x), P(\cdot,\cdot), Q(\cdot,\cdot)$ are continuous and continuously differentiable in their domain, and conditions (\ref{eq22}), (\ref{eq24}), $\gamma_1\geq0$, $k>0$ are fulfilled, then formulated problem has a unique solution represented as follows
\[
\begin{array}{l}
u(x,t)=\theta(t)\frac{2}{l}\Gamma[\alpha+\mu(1-\alpha)]\int\limits_0^l\left[F_2(\xi)+\int\limits_0^l F_2(\eta)R(\xi,\eta)d\eta\right]t^{(1-\mu)(\alpha-1)}\sum\limits_{n=1}^\infty E_{\alpha,\alpha+\mu(1-\alpha)}\left(-k\lambda_n^2t^\alpha\right)\times\\\sin(\lambda_nx)\sin(\lambda_n\xi)d\xi+\frac{\theta(-t)}{2}\left[F_2(x+t)+F_2(x-t)+\int\limits_0^l F_2(\xi)[R(x+t,\xi)+R(x-t,\xi)]d\xi+\right.\\
\left.\int\limits_{x-t}^{x+t}\left[F_2'(z)-2\psi'(z)+\int\limits_0^l[F_2'(\eta)-2\psi'(\eta)]R(z,\eta)d\eta\right]dz\right],
\end{array}
\]
where $\theta(t)=1$ for $t\geq 0$ and $\theta(t)=0$ for $t<0$.

\textsc{Remark.} Functions $P(x,t)=\sin x\cos t$ and $Q(x,t)=e^{-x}\left(1+e^{-t}\right)$ satisfy all conditions imposed in Theorem if $l\leq \pi/2$.

%

\end{document}